\documentclass{article}
\usepackage{amssymb,latexsym,amsmath,amsthm}
\usepackage{fancyhdr}
\usepackage{color}
\usepackage{mathrsfs}

\newtheorem{thm}{Theorem}[section]
\newtheorem{prop}[thm]{Proposition}
\newtheorem{lemma}[thm]{Lemma}
\newtheorem{cor}[thm]{Corollary}

\newcommand{\ep}{\varepsilon}

\newcommand{\con}{\equiv}

\newcommand{\ndiv}{\nmid}
\newcommand{\modd}[1]{\; ( \text{mod} \; #1)}
\newcommand{\bstack}[2]{#1 \atop #2}

\newcommand{\maps}{\rightarrow}

\newcommand{\al}{\alpha}

\newcommand{\del}{\delta}

\newcommand{\om}{\omega}

\newcommand{\sig}{\sigma}
\newcommand{\lam}{\lambda}
\newcommand{\Lam}{\Lambda}

\newcommand{\Pcal}{\mathcal{P}}

\newcommand{\Q}{\mathbb{Q}}

\newcommand{\Z}{\mathbb{Z}}

\newcommand{\beq}{\begin{equation}}
\newcommand{\eeq}{\end{equation}}
\numberwithin{equation}{section}

	\newcommand{\xtra}[1]{}
\newcommand{\Cl}{\mathrm{Cl}}
\newcommand{\pfk}{\mathfrak{p}}

\begin{document}

\title{Averages and moments associated to class numbers of imaginary
  quadratic fields}  
\author{D. R. Heath-Brown\footnote{Mathematical Institute, Radcliffe
    Observatory Quarter, Woodstock Road, Oxford 
OX2~6GG, {\tt rhb@maths.ox.ac.uk}} \\
	  L. B. Pierce\footnote{Department of Mathematics, Duke
            University, Durham NC 27708, {\tt pierce@math.duke.edu}
            } }
  \date{}
\maketitle  
\begin{abstract}
For any odd prime $\ell$, let $h_\ell(-d)$ denote the $\ell$-part of the class number of the imaginary quadratic field $\Q(\sqrt{-d})$. Nontrivial pointwise upper bounds are known only for $\ell =3$; nontrivial upper bounds for averages of $h_\ell(-d)$ have previously been known only for $\ell =3,5$. In this paper we prove nontrivial upper bounds for the average of $h_\ell(-d)$ for all primes $\ell \geq 7$, as well as nontrivial upper bounds for certain higher moments for all primes $\ell \geq 3$.
\end{abstract}
\noindent
%%The subject numbers below were submitted to journal but were excluded for arXiv due to formatting issues
%\nonumfootnote{\emph{2010 Mathematics Subject Classification:} 11R29, 11D45}
%\nonumfootnote{\emph{Keywords:} class numbers, Cohen-Lenstra heuristics}
%%11R29: Class numbers, class groups, discriminants
%%11D45: Counting solutions of Diophantine equations

\section{Introduction}
Fix an imaginary quadratic field $\Q(\sqrt{-d})$ with square-free
$-d<0$, and let $\Cl(-d)$ be the corresponding class group. The size
of the class group, denoted $h(-d)$, is the class number of
$\Q(\sqrt{-d})$, a fundamental invariant that appears widely in number
theory. The divisibility properties of class numbers of quadratic
fields are subject to the conjectures known as the Cohen-Lenstra
heuristics \cite{CohLen84}, which despite significant attention remain open in most
cases. For any prime $\ell \geq 2$, let $h_\ell(-d)$ denote the $\ell$-part of
the class number, that is the number of ideal classes in the class
group $\Cl(-d)$ whose $\ell$-th power is the principal ideal class. One
may obtain a trivial pointwise upper bound for $h_\ell(-d)$ by noting
that  
\[ h_\ell(-d) \leq h(-d) \ll d^{1/2+ \ep}.
\]
It is conjectured that 
\beq\label{h_conj}
h_\ell(-d) \ll d^{\ep} 
\eeq
for all $d$ and any $\ep>0$. (Throughout, we will use the convention that all implied constants may depend upon $\ell$ and $\ep$.)

This conjecture (and a more general version for $\ell$-torsion in class groups of number fields of any degree) is motivated by the Cohen-Lenstra heuristics \cite{CohLen84}, by counting elliptic curves with  fixed conductor  \cite{BruSil96}, by counting number fields of fixed degree and discriminant \cite{Duk98}, and by questions on equidistribution of CM-points on Shimura varieties \cite{Zha05}.
For $\ell =2$, the conjecture (\ref{h_conj})  is known by the genus
theory of Gauss. 
For $\ell =3$ the currently best known upper bound is due
to Ellenberg and Venkatesh \cite{EllVen07}: 
 \beq\label{EllVen}
 h_3(-d) \ll d^{1/3 +\ep}.
 \eeq
  For  primes $\ell \geq 5$, no nontrivial upper bound for $h_\ell(-d)$ is
  known to hold for all $d$.  
 
 One may also consider averages
\[ \sum_{0 < d< X} h_\ell(-d).\]
In the case $\ell =3$, Davenport and Heilbronn \cite{DavHei71}  established that
\beq\label{DavHeil}
 \sum_{0 < d< X} h_3(-d)  \sim 2 \sum_{0<d<X} 1,
 \eeq
 as $X \maps \infty$,
 in which both sums are restricted to fundamental discriminants.
This asymptotic has recently been refined further to include secondary
main terms (see Bhargava, Shankar and Tsimerman \cite{BST}, Taniguchi
and Thorne \cite{TanTho13}, and Hough \cite{Hou13}), but for the 
purposes of this paper it is sufficient that (\ref{DavHeil}) provides
an upper bound: 
\beq\label{DavHeil_bound}
 \sum_{0<d<X} h_3(-d) \ll X.
\eeq
  For $\ell = 5$, the best known upper bound for the average is due to
  Soundararajan \cite{Sou00} (also proved by Hough \cite{Hou13}): 
\beq\label{Sound5}
 \sum_{0<d<X} h_5(-d) \ll X^{5/4+\ep}.
 \eeq
\xtra{
Extra note on Sounds' method: Sound shows in his paper that $h_\ell(-d)
-1 \leq 2D^+$ where roughly speaking  
\[ D^+ = \# \{ lm^\ell = n^2 l^2 + t^2 d, l | d, (m,2dn)=1, m< C^+
\sqrt{d(\al l)^{-1}} \}.\] 
Supposing we only consider prime $d$ then 
\[ D^+ \lesssim\# \{ m^\ell = n^2 + t^2 d,  m< C^+ \sqrt{d} \}.\]
Thus we are now counting solutions to 
\[ m^\ell = n^2 +t^2 d\]
with $m= O(\sqrt{d})$, $n \leq N= O(d^{\ell/4})$ and $t =
O((m^\ell/d)^{1/2}) = O(d^{(\ell-2)/4})$. We count the number of solutions
to this equation on average over $X \leq d < 2X$ by counting
solutions to the congruence $m^\ell \con n^2 \modd{t^2}$, leading to an
upper bound for $D^+$ of the form 
\[ \sum_{m} \left( \frac{N}{t^2} +1 \right) N^\ep (\# t) \ll X^{1/2}
\left( \frac{X^{\ell/4}}{X^{(\ell-2)/2}} +1 \right) X^\ep (X^{(\ell-2)/4})  
	\ll X^{\ell/4+\ep} + X^{1+\ep} \ll X^{\ell/4+\ep}\]
	for $\ell \geq 5$. This is nontrivial if $\ell =5$ and trivial for
	all $\ell \geq 6$. Hence Sound only states a nontrivial bound for
	the case $\ell =5$. 
}
For primes $\ell \geq 7$, the literature appears to contain no bound better than the 
 trivial estimate
\[ \sum_{0<d<X} h_\ell(-d) \ll X^{3/2+\ep}.\]
However Soundararajan noted in \cite{Sou00} that he has shown for any
prime $\ell \geq 3$ that  
\beq\label{h_small} h_\ell(-d) \ll d^{\frac{1}{2} - \frac{1}{2\ell}+\ep}
\eeq
 for all but one square-free discriminant $d$ in any dyadic range
 $[X,2X)$. Summing over $O(\log X)$ dyadic ranges implies the
 nontrivial average bound 
\beq\label{h_avgt}
 \sum_{0<d<X} h_\ell(-d) \ll X^{\frac{3}{2} - \frac{1}{2\ell} +\ep} \eeq
for any $\ell \geq 3$.
While this
is superseded by (\ref{DavHeil_bound}) and (\ref{Sound5}) for $\ell =3$
and $5$,  no improvement
has been given hitherto for larger values of $\ell$.

One can further consider the second moment; motivated by the
conjecture (\ref{h_conj}) for the pointwise upper bound for $h_\ell(-d)$,
one would expect that  
\[ \sum_{0<d<X} h_\ell(-d)^2 \ll X^{1+\ep}.\]
For $\ell =3$ and $5$, one may bound the second moment by applying the
best known pointwise upper bound  (respectively (\ref{EllVen}) and
(\ref{h_small}))  to one factor $h_\ell(-d)$, and then applying the best
known average upper bound  to the remaining sum (respectively
(\ref{DavHeil_bound}) and (\ref{Sound5})). For $\ell \geq 7$, it is
advantageous to apply Soundararajan's result (\ref{h_small}) to both
factors of $h_\ell(-d)$. This approach results in the following upper
bounds for the second moment: 
\beq\label{moment_1}
  \sum_{0<d<X} h_\ell(-d)^2 \ll \begin{cases}
	X^{\frac{4}{3}+\ep} & \ell =3 \\
	X^{\frac{33}{20} + \ep} & \ell =5 \\
	X^{2 - \frac{1}{\ell} + \ep} & \ell \geq 7, \; \text{prime.}
	\end{cases}
	\eeq
More generally, for any real number $k \geq 1$, known results lead to
bounds for the $k$th moment of the form 
\[
  \sum_{0<d<X} h_\ell(-d)^k \ll \begin{cases}
	X^{1 + \frac{(k-1)}{3}+\ep} & \ell =3 \\
	X^{\frac{5}{4} + (k-1)(\frac{2}{5}) + \ep} + X^{\frac{k}{2}
          +\ep} & \ell =5 \\ 
	X^{1+k( \frac{\ell-1}{2\ell}) + \ep} + X^{\frac{k}{2}+\ep} & \ell \geq
        7, \; \text{prime.} 
	\end{cases}
	\]

\subsection{Statement of the Theorems}
The purpose of this paper is to improve on these bounds for the
averages and moments of $h_\ell(-d)$ for $d$ square-free and $\ell$ an odd prime.
(For the rest of this paper the notations $d$ and $\ell$ are reserved for
square-free integers and odd primes respectively.)

\begin{thm}\label{thm_avg}
For each prime $\ell \geq 5$,
\[  \sum_{0<d<X} h_\ell(-d) \ll X^{\frac{3}{2} - \frac{3}{2\ell+2} +\ep},\]
for any $\ep>0$.
\end{thm}
This recaptures Soundararajan's result (\ref{Sound5}) for $\ell =5$ and
improves on the bound (\ref{h_avgt}) for all primes $\ell \geq 7$.  
(Since Davenport and Heilbronn's result (\ref{DavHeil}) is best
possible, our work provides no new information for the average of
$h_3(-d)$.) 

We also consider higher moments. First we consider the moments
of $h_3(-d)$, for which our main result is the following: 
\begin{thm}\label{thm_square}
\[  \sum_{0<d<X} h_3(-d)^4 \ll X^{\frac{11}{6}+\ep} \qquad \text{for any $\ep>0$}. \] 
\end{thm}
It may be surprising to see the 4th moment here, but it turns out to
give the best results of its type, as we shall see.

By the reflection principle of Scholz \cite{Sch32}, $\log_3 h_3(-d)$ and
$\log_3 h_3(+3d)$ differ by at most one.  
Thus the corresponding bound for the $3$-part of the class number of
real quadratic fields follows as a corollary, making an identical
improvement over previously known bounds as in the imaginary case: 
\begin{cor}\label{cor_pos}
\[  \sum_{0<d<X} h_3(d)^4 \ll X^{\frac{11}{6}+\ep} \qquad \text{for any $\ep>0$}.\]
\end{cor}

Nontrivial bounds for other moments are also an immediate corollary.
For $1\le k<4$ one merely uses H\"older's inequality in conjunction
with (\ref{DavHeil_bound}), while for $k>4$ one just applies (\ref{EllVen}):
\begin{cor}\label{cor_higher}
For all real $k\in[1,4]$, and for any $\ep>0$,
\begin{eqnarray*}
\sum_{0<d<X} h_3(-d)^k & \ll & X^{(5k+13)/18+\ep} \\
\sum_{0<d<X} h_3(d)^k & \ll &X^{(5k+13)/18+\ep}.\\
\end{eqnarray*}
For all real $k\ge 4$, and for any $\ep>0$,
\begin{eqnarray*}
\sum_{0<d<X} h_3(-d)^k & \ll & X^{(2k+3)/6+\ep} \\
\sum_{0<d<X} h_3(d)^k & \ll &X^{(2k+3)/6+\ep}.\\
\end{eqnarray*}
In particular, for any $\ep>0$,
\[\sum_{0<d<X} h_3(-d)^2  \ll  X^{23/18+\ep}. \]
\end{cor}
This final bound improves on (\ref{moment_1}); 
we note that $23/18 = 1.2777...$.

We next consider higher moments for $h_\ell(-d)$ for primes $\ell \geq 5$.
Theorem \ref{thm_avg} combined with (\ref{h_small}) implies that for
any real $k \geq 1$, 
\[
 \sum_{0<d<X} h_\ell(-d)^k \ll X^{\frac{3}{2} - \frac{3}{2\ell+2} +
   (k-1)(\frac{1}{2} - \frac{1}{2\ell}) + \ep} + X^{\frac{k}{2}+\ep}, 
 \]
where the last term arises from the possible exceptions to (\ref{h_small}). 
For purposes of comparison, we rewrite this as
\[  \sum_{0<d<X} h_\ell(-d)^k \ll X^{ 1 + k \left( \frac{\ell-1}{2\ell} \right) -
    \frac{2\ell-1}{2\ell(\ell+1)} + \ep} + X^{\frac{k}{2}+\ep}. 
  \]
We will improve on this for all real $1 < k < (2\ell^2+1)/(\ell+1)$:
\begin{thm}\label{thm_k_real}
For any prime $\ell \geq 5$, all real $k\ge 1$, and any $\ep>0$,
\[ \sum_{0<d<X} h_\ell(-d)^k \ll 
	\begin{cases}
	X^{1+ k\left( \frac{\ell-2}{2\ell+2} \right) + \ep} & \text{if $1
          \leq k \leq \frac{\ell^2-1}{2\ell-1}$},\\ 
	X^{1 +k\left( \frac{\ell-1}{2\ell} \right) - \left( \frac{\ell-1}{2\ell}
          \right)  + \ep} & \text{ if $\frac{\ell^2-1}{2\ell-1} \leq k \leq
          \ell+1$},\\
X^{\frac{k}{2}+\ep}& \text{ if $k\ge \ell+1$}.
	\end{cases}\]
\end{thm}
In particular, we single out the consequence of Theorem
\ref{thm_k_real} for the second moment (noting that $k=2$ lies in the
first case of the theorem for $\ell \geq 5$): 
\begin{cor}
For any prime $\ell \geq 5$, for any $\ep>0$,
\[ \sum_{0<d<X} h_\ell(-d)^2 \ll X^{2 - \frac{3}{\ell+1} + \ep}.\]
\end{cor}
This improves on (\ref{moment_1}) in every case.
Theorem \ref{thm_avg} may of
course be deduced from the above corollary via the Cauchy-Schwarz
inequality.  However we have stated and proved Theorem \ref{thm_avg} separately
since it is, in effect, used in the proof of Theorem  \ref{thm_k_real}.

Our approach is to develop an unconditional upper bound for $h_\ell(-d)$
that holds for almost all $d$, by using the relation between $h_\ell(-d)$
and small split primes in
$\Q(\sqrt{-d})$. The original observation of this relation is 
credited to Soundararajan (and to Michel in a related context) in work of Helfgott and Venkatesh \cite{HelVen06} and Ellenberg and Venkatesh \cite{EllVen07},
and has been used in \cite{HelVen06}, for example, to prove a bound for
$h_3(-d)$ for all $d$, conditional on GRH. Here we
prove an unconditional version, at the cost that it only holds for
``almost all'' $d$. To treat higher moments, we combine this with
upper bounds for the number of simultaneous representations of
integers by certain polynomials; this counting problem is similar to
computations performed in \cite{Sou00} and \cite{HB07a}. Finally, we
remark that the methods of Section \ref{sec_higher} may also be 
applied to prove upper bounds for mixed averages of the form  
\[ \sum_{0<d<X} h_\ell(-d) h_{\ell'}(-d) \]
for distinct odd primes $\ell,\ell'$; we leave the details to the interested reader.

We reiterate that throughout this paper we consider sums over $0<d<X$ to be
restricted to square-free integers, and $\ell$ represents an odd
prime. We will frequently combine factors of size $X^\ep$ for various
$\ep$; in all cases $\ep$ may be taken to be an arbitrarily small real
number, so we re-define it wherever appropriate so that the total
factor remains represented by $X^\ep$. We also use the notation $A\ll
B$ to indicate that there is a constant $c$, possibly depending on
certain allowable parameters such as $\ell$ or $\ep$, 
such that $|A| \leq c|B|$,
and similarly for $A \gg B$. 

\section{An unconditional pointwise upper bound}\label{sec_ptwise}
Our starting point is the following unconditional pointwise upper
bound for $h_\ell(-d)$: 
\begin{prop}\label{prop_hg_bound}
Fix any prime $\ell \geq 3$ and real parameters $\tfrac14 X^{\frac{1}{2\ell}} \le Z
\le  X$. 
There
exists a small exceptional set $E(Z;X) \subset [X,2X)$ such that for
all square-free $d \in [X,2X) \setminus E(Z;X)$, 
\[ h_\ell(-d) \ll X^\ep \left\{ d^{1/2} Z^{-1} + d^{1/2} Z^{-2}S_\ell(d;Z) \right\}, 
 \]
for any $\ep>0$, where $S_\ell(d;Z)$ is the cardinality of the 
set of pairs of primes $p,p'$ satisfying
\[ Z \leq p \neq p'  < 2Z\]
for which there exist $u,v \in \Z \setminus \{0\}$ with $(v,pp')=1$ such that 
\[ 4(pp')^\ell = u^2+dv^2. \]
Moreover,  the exceptional set satisfies
\beq\label{E_small}
 \# E(Z;X) \ll X^{\ep'}
 \eeq
 for any $\ep'>0$.
\end{prop}

\begin{cor}\label{newcor}
Fix any $\ep'>0$. For all $d\in [X,2X)$ apart from at most $O(X^{\ep'})$ exceptions,
\[ h_\ell(-d)\ll d^{\frac{1}{2}-\frac{1}{2\ell}+\ep}\]
for any $\ep>0$.
\end{cor}
 This corollary, which we will prove at the end of Section \ref{sec_ptwise}, gives a weak form of Soundararajan's result concerning
the bound (\ref{h_small}).

It is clear from Proposition \ref{prop_hg_bound}
 that an understanding of $S_\ell(d;Z)$, both
in terms of its average over $d$ and its second moment, will yield
corresponding information for $h_\ell(-d)$. Our two main technical
results are for the average and second moment of  $S_\ell(d;Z)$: 
\begin{prop}\label{prop_Sg_avg}
For any prime $\ell \geq 3$ and $X^{\frac{1}{2\ell}} \le Z \le X$, 
\[ \sum_{X \leq d < 2X}S_\ell(d;Z)
	\ll  X^\ep \{ Z^2X^{1/2} +Z^{\ell+2}X^{-1/2}\}
	\]
	for any $\ep>0$.
\end{prop}

\begin{prop}\label{prop_Sg_moment}
For $\ell =3$ and $X^{\frac{1}{6}}  \le Z \le X$, 
\[ \sum_{X \leq d < 2X} S_3(d;Z)^2 \ll 
X^\ep \{ Z^2X^{1/2} + Z^{12}X^{-3/2} \}\]
for any $\ep>0$.
For any prime $\ell \geq 5$ and 
$X^{\frac{1}{2\ell}}  \le Z \le X$, 
\[\sum_{X \leq d < 2X} S_\ell(d;Z)^2 \ll
X^\ep \{ Z^2X^{1/2} + Z^{2\ell+4}X^{-1} \} \]
for any $\ep>0$.
\end{prop}
We include the case $\ell \geq 5$ in Proposition \ref{prop_Sg_moment} as it requires little extra effort, but we will not make use of it: while it does result in a nontrivial upper bound for the second moment of $h_\ell(-d)$, a stronger result may be obtained by applying Proposition \ref{prop_Sg_avg} directly.

 In the remainder of this section, we prove Proposition
 \ref{prop_hg_bound} and its corollary. 
We prove Propositions \ref{prop_Sg_avg} and
 \ref{prop_Sg_moment} in Sections \ref{sec_avg_S} and \ref{sec_T1T2},
 respectively. 
Finally, in Sections \ref{sec_avg_h} and \ref{sec_higher} we record the
consequences of these results for averages and moments of $h_\ell(-d)$.

\subsection{Proof of Proposition \ref{prop_hg_bound}}
Fix a prime $\ell \geq 3$ and a square-free integer $X \leq d < 2X$. Let
$H=\Cl(-d)$ be the class group of $\Q(\sqrt{-d})$, with class number
$h(-d) = \#\Cl(-d)$.  
Let $H_\ell$ denote the maximal elementary abelian $\ell$-group
in $H$, with $h_\ell(-d) = \# H_\ell$. Since 
\beq\label{hHH}
 \# H/H_\ell = \frac{h(-d)}{h_\ell(-d)} ,
 \eeq
in order to show that $h_\ell(-d)$ is small it suffices to show that
there are many cosets of $H_\ell$ in $H$. Let $\chi_d(\cdot)$ denote the
quadratic character associated to $\Q(\sqrt{-d})$. Picking a prime $p
\ndiv 2d$ 
such that $\chi_d(p)=1$, it follows that $p$ splits in
$\Q(\sqrt{-d})$ as $\pfk \pfk^\sig$, say, where $\sig$ is 
the  non-trivial Galois
automorphism of $\Q(\sqrt{-d})$. Suppose that two distinct primes $p,
p'$ split in this manner as $\pfk \pfk^\sig$ and $\pfk' \pfk'^\sig$
respectively, and suppose that $\pfk$ and $\pfk'$ represent the same class
in $H/H_\ell$, so that $\pfk H_\ell = \pfk' H_\ell$. It follows that
$\pfk^{-1}\pfk' \in H_\ell$, so that $(\pfk^{-1} \pfk')^\ell$ is a principal
ideal. Thus  
$(\pfk^\sig \pfk')^\ell$ is also a principal ideal, say  
\beq\label{ideal_id}
 (\pfk^\sig \pfk')^\ell= \left( \frac{u+ v\sqrt{-d}}{2}\right),
 \eeq
for some $u,v \in \Z$. 
Hence taking norms, it follows that 
\beq\label{ell_eqn}
 4 (pp')^\ell = u^2 +dv^2 .
 \eeq
Note that we may require that $\gcd (v,pp')=1$ (and in particular that
$v \neq 0$). For supposing that $p|v$, say, then by (\ref{ell_eqn}) we
see that also $p|u$ so that 
$p | ( \frac{u+ v\sqrt{-d}}{2})$. Hence 
$\pfk | ( \frac{u+ v\sqrt{-d}}{2})$, which
by (\ref{ideal_id}) implies that $ \pfk | 
(\pfk^\sig\pfk')^\ell$. Since $p$ is unramified this would then imply
that $\pfk | \pfk'$, which 
contradicts the fact that $p \neq p'$. A similar argument shows that
we may require that $u \neq 0$.  

We will show that for all but a small number of ``exceptional'' $d$,
there are many primes $p, p'$ that split in this manner, while also
showing there can only be few solutions $(u,v)$ to (\ref{ell_eqn}) with
$\gcd(v,pp')=1$ and $u,v$ in an appropriate range. This forces there
to be many distinct cosets of $H_\ell$ in $H$, and provides an upper
bound for $h_\ell(-d)$, as long as $d$ is not exceptional.

We first fix $X \leq d<2X$ and count the number of primes $p$ that split
appropriately, with  
\[ Z \leq p< 2Z\]
 for some parameter $Z$ with  $\tfrac14 X^{1/2\ell}  \le Z \le X$ (to be chosen
 precisely in applications). We see that  
\[ \#\{ Z  \leq p < 2Z: \chi_d(p) = 1 \} 
= \frac{1}{2} \sum_{Z \leq p < 2Z} (1 + \chi_d(p)) + O(\om(d)),\]
where the last term reflects the contribution of the primes that divide $d$, and contributes no more than $O(\log X) = O(\log Z)$. We now separate the two terms within the sum over $p$ and apply the prime number theorem, obtaining
\[ \#\{ Z  \leq p < 2Z: \chi_d(p) = 1 \} = \frac{1}{2} Z(\log Z)^{-1} +\frac{1}{2}M(d;Z)+ O(Z(\log Z)^{-2}),\]
say, where 
\[ M(d;Z) = \sum_{Z \leq p < 2Z} \chi_d(p).\]
Thus the number of split primes in this range is $\gg Z(\log Z)^{-1}$,
unless we have $|M(d;Z)| \ge \tfrac14 Z (\log Z)^{-1}$; 
we will show this exceptional
scenario can occur for only a small number of $d$.  

Given a character $\chi$, set
\[ V(\chi) = (\sum_{Z \leq  p < 2Z} \chi(p))^{4\ell}.\]
Upon unfolding the product, we see that this is a character sum of the form 
\[ \sum_{Z^{4\ell}  \leq n < (2Z)^{4\ell}} a_n \chi(n)\]
 for some coefficients $|a_n| \ll d(n)^{4\ell} \ll Z^\ep$.
 Now we note that with the particular choice $\chi = \chi_d$,
\beq\label{SV}
 \sum_{X \leq d< 2X} |M(d;Z)|^{8\ell}  = \sum_{X \leq d< 2X} |V(\chi_d)|^2.
 \eeq
 By positivity, we can enlarge the sum on the right hand side of 
 (\ref{SV}) to include all primitive characters modulo $d$ and apply
 the large sieve (see for example \cite[Thm. 4, Ch. 27]{Dav}), to obtain 
\begin{eqnarray}
  \sum_{X \leq d <2 X} |M(d;Z)|^{8\ell} & \leq& \sum_{X \leq d <2X} \;
  \sideset{}{^*}\sum_{\chi \modd{d}} |V(\chi)|^2 \nonumber \\ 
  & \ll &(X^2 + Z^{4\ell})(\sum_{Z^{4\ell} \leq n < (2Z)^{4\ell}} |a_n|^2) \nonumber \\
	&\ll& Z^{4\ell+2\ep}(X^2+Z^{4\ell}) \ll Z^{8\ell + 2\ep}, \label{largesieve}
	\end{eqnarray}
	since $X^{1/2\ell} \ll Z$ by assumption.
Let $E(Z;X)$ denote the exceptional set, 
\beq\label{exp_set}
 E(Z;X) = \{ X \leq d < 2X: |M(d;Z)|\ge\tfrac14 Z(\log Z)^{-1}\}.
 \eeq
 Then we may conclude from (\ref{largesieve}) that the exceptional set is small:
\[ \#E(Z;X) \ll X^{\ep},
\]
for any $\ep>0$.

We now fix a $d$ with $X\leq d <2X$ such that $d \not\in E(X,Z)$; the above
argument shows that there are $\gg Z(\log Z)^{-1}$ split primes for
this $d$. In particular, summing over all cosets of $H_\ell$ in $H$ shows
that for this $d$, 
\begin{eqnarray*}
\lefteqn{ \sum_{C \in H/H_\ell} \#\{ Z \leq p < 2Z : \chi_d(p)=1, p=\pfk\pfk^\sig,
 \pfk \in C\}} \\ 
& = & \#\{ Z \leq p < 2Z : \chi_d(p)=1\} \gg Z(\log Z)^{-1}.
 \end{eqnarray*}
On the other hand, applying the Cauchy-Schwarz inequality to the left
hand side shows that
\beq\label{HH_upper}
 (\# H/H_\ell)^{1/2}  \left( S_\ell^{(1)}(d;Z) \right)^{1/2}\gg Z(\log Z)^{-1},
 \eeq
where we define
\[ S^{(1)}_\ell(d;Z) = \sum_{C \in H/H_\ell}\#\{ Z \leq p < 2Z: \chi_d(p)=1,
p = \pfk \pfk^\sig, \pfk \in C\}^2 .\] 
By the above discussion, we know that
\beq\label{S1_set}
 S^{(1)}_\ell(d;Z) \ll \#\{ Z \leq p, p' < 2Z: 4(pp')^\ell = u^2+dv^2 \;
 \text{for some $u,v \in \Z$} \}, 
 \eeq
 where in the case that $p \neq p'$ we may impose the additional
 conditions that $u,v \neq0$ and $(v,pp')=1$. 
Combining (\ref{HH_upper}) and (\ref{hHH}), we may conclude that
\[h_\ell(-d) \ll d^{1/2+\ep} Z^{-2} (\log Z)^{2} S^{(1)}_\ell(d;Z),
\]
still under the assumption that $d$ is not exceptional.
Finally, we write 
\[ S^{(1)}_\ell(d;Z) = S_\ell^{(0)}(d;Z) + S_\ell(d;Z),\]
 where $S_\ell^{(0)}(d;Z) $ is the contribution to the set (\ref{S1_set})
 from pairs $p=p'$ and $S_\ell(d;Z)$ is the contribution from pairs
 $p \neq p'$. Trivially, $S_\ell^{(0)}(d;Z)  \ll Z$, and we see that
 Proposition \ref{prop_hg_bound} holds. 

To deduce the corollary we take $Z=\tfrac14 X^{1/2\ell}$, and note that
any pairs of primes $p,p'$ counted by $S_\ell(d;Z)$ would satisfy
\[X\le d\le u^2+dv^2=4(pp')^\ell\le 4(4Z^2)^\ell=4^{1-\ell}X<X.\]
Thus $S_\ell(d;Z)$ must vanish, so that $h_\ell(-d)\ll X^{\ep}d^{1/2}Z^{-1}
\ll d^{1/2-1/(2\ell)+\ep}$ unless $d$ lies in $E(Z;X)$.  The result then follows.

\section{Proof of Proposition \ref{prop_Sg_avg}}\label{sec_avg_S}
Define the parameters
\beq\label{param}
 W= Z^2,\qquad U = 2^{\ell+1} Z^\ell, \qquad V = 2^{\ell+1}Z^\ell X^{-1/2}.
 \eeq
 Note that $V \ge 2$ as long as 
 \[ Z \ge X^{\frac{1}{2\ell}},
 \]
  which we henceforward assume.
 Note also that up to a constant factor (accounting for changing signs of $u,v$), we may express $S_\ell(d;Z)$ as the quantity
\[  \#\{ Z \leq p \neq p' < 2Z: 4(pp')^\ell = u^2+dv^2 \;
 \text{for some $u,v \geq 1$ with $(v,pp')=1$} \}.\]
Furthermore, for any $X \leq d< 2X$, any triple $w=pp'$, $u,v$ considered in the set above
 certainly satisfies $W \leq w < 4W$, $1 \leq u \leq U$, $1 \leq v \leq V$. 

We wish to bound $S_\ell(d;Z)$ on average over $d$; for this we note that 
\begin{eqnarray*}
 \sum_{\bstack{X\leq d < 2X}{d \not\in E(Z;X)}}  S_\ell(d;Z) & \ll & \#\{W \leq w < 4W, 1 \leq u \leq U, 1 \leq v \leq V:
 \gcd(v,w)=1, \nonumber \\ 
&& \hspace{1cm} v^2 | (4w^\ell - u^2), (4w^\ell - u^2)/v^2 \in [X,2X) \}.
	  \end{eqnarray*}
It is convenient to work with dyadic ranges; thus for any parameter $1 \leq V_0 \leq V/2$, define
\begin{eqnarray*}
 N(Z,X;V_0) &=& \#\{W \leq w < 4W, 1 \leq u \leq U, V_0\leq v < 2 V_0:
 \gcd(v,w)=1, \nonumber \\ 
&& \hspace{1cm} v^2 | (4w^\ell - u^2), (4w^\ell - u^2)/v^2 \in [X,2X) \}.
	  \end{eqnarray*}
Then certainly
	  \[ \sum_{\bstack{X\leq d < 2X}{d \not\in E(Z;X)}}  S_\ell(d;Z)  \ll 
	  	\sum_{0 \leq j \leq \log_2 (V)-1} N(Z,X;2^j)  = 
	  \sum_{\bstack{V_0 \leq V/2}{\text{dyadic}}} N(Z,X;V_0).
	  \]
	  
We  turn to bounding an individual term $N(Z,X;V_0)$.
We first fix $w$ and $v$ and let 
\[ M(w;v) = \# \{ u  \modd{v^2} : u^2 \con 4w^\ell\modd{v^2}\}.\]
\begin{lemma}\label{lemma_M}
For any coprime $w$ and $v$, 
\beq\label{M_total}
M(w;v) \leq 2^{\omega(v)+1} \ll v^\ep,
\eeq
where $\omega(v)$ denotes the number of distinct prime divisors of $v$.
\end{lemma}
This is proved in a standard fashion.
  Writing $v=q_1^{r_1} \cdots q_s^{r_s}$ in its  prime decomposition,
  it suffices by the Chinese Remainder Theorem to count
  $M(w;q_i^{r_i})$ for each $q_i $. 
Since $(w,v)=1$ we may assume that $(w,q_i)=1$; we
  also assume for the moment that $q_i$ is odd. Then $M(w;q_i^{r_i})$
  will be nonzero only if $w$ is a quadratic residue modulo $q_i$, in
  which case $u$ can lie in at most 2 residue classes modulo $q_i$;
  since $q_i$ is odd, each solution modulo $q_i$ lifts uniquely to a
  solution modulo $q_i^{2r_i}$. Thus we see that in this case 
\[M(w;q_i^{r_i})\leq 2.
\]
 If $q_i=2$ then the relevant congruence has solutions only if $2|u$,
 in which case we may equivalently count solutions to $
 (u/2)^2 \con w^\ell \modd{q_i^{2r_i-2}}$.
However if $n$ is odd, a congruence $x^2\con n\modd{2^r}$ has at most 4
solutions. We may therefore conclude that $M(w;q_i^{r_i})\le 4$,
 thus proving (\ref{M_total}).

Applying Lemma \ref{lemma_M} directly to count solutions $u \leq U$ to
$u^2 \con 4 w^\ell\modd{v^2}$ would lead to the upper bound  
\beq\label{NV_weak}
 N(Z,X; V_0) \ll WV_0^{1+\ep}(UV_0^{-2} +1).
 \eeq
But then summing over all dyadic ranges with $1 \leq V_0 \leq V/2$
would not allow us to take advantage of the decay with respect to
$V_0$ in (\ref{NV_weak}). Thus we return to the definition of
$N(Z,X;V_0)$ and utilize the additional piece of information that  
\[ X \leq \frac{4w^\ell-u^2}{v^2} < 2X,\]
which we re-write as 
\beq\label{v_interval}
v^2X \leq 4w^\ell - u^2 <2 v^2 X.
\eeq
 We will conclude from this that $u$ must lie within a short 
interval around $2w^{\ell/2}$; precisely, we write
 \[\left(\frac{u}{2w^{\ell/2}}\right)^2=1+E,\]
in which (\ref{v_interval}) shows that
 \[|E|\leq\frac{2Xv^2}{4w^\ell}\le
 \frac{8XV_0^2}{4W^\ell}=\frac{2XV_0^2}{Z^{2\ell}}=2^{2\ell+3}\frac{V_0^2}{V^2}.\]
 Thus $E\ll 1$ whence $\sqrt{1+E}=1+O(E)$.  It follows that
 \[u=2w^{\ell/2}+O(w^{\ell/2}E)=2w^{\ell/2}+O(W^{\ell/2}V_0^2V^{-2}).\] 
\xtra{Old argument: this implies that for $u,v,w$ counted by $N(Z,X;V_0)$, 
\beq\label{uwVX}
 u^2 = w^\ell + E_0,
 \eeq
 where $|E_0| < 8V_0^2X$. 
Recall that $Z^{2\ell} \leq w^\ell \leq 4 Z^{2\ell}$ so that $V^2X \leq 2^{-2\ell}
w^\ell$.  After multiplying through by $V_0^2 V^{-2}$ we see that $V_0^2
X \leq 2^{-2\ell} w^\ell V_0^2 V^{-2}$.  We can now rewrite (\ref{uwVX}) as 
 \[ u^2 = w^\ell(1 + E_1) ,\]
 where 
 \[ |E_1| < 2^{-2\ell+3}V_0^2V^{-2}.\]
Since $V_0/V \leq 1/2$ we may conclude that $|E_1| <1$ and so certainly
\[ u = w^{\ell/2}(1 + E_1).\]
}

Thus for each fixed $w,v$, in order to be counted by $N(Z,X;V_0)$, $u$
must lie in an interval $I_w$ around $2w^{\ell/2}$ of length $O(W^{\ell/2}V_0^2 V^{-2})$. We
apply this information along with the bound (\ref{M_total}) to
conclude that for each fixed $w,v$ considered in $N(Z,X;V_0)$,  
\[ \#\{ u \in I_w : u^2 \con 4w^\ell \modd{v^2} \}
	\ll  V_0^\ep \left(\frac{W^{\ell/2} V_0^{2} V^{-2}}{V_0^2} + 1 \right)
	=  V_0^\ep (W^{\ell/2} V^{-2} + 1).\]
As a consequence, 
\begin{eqnarray}
N(Z,X;V_0)& \ll &\sum_{\bstack{W \leq w < 4W, V_0 \leq v <2
    V_0}{(v,w)=1}}  \#\{ u \in I_w : u^2 \con 4w^\ell \modd{v^2} \}
\nonumber\\ 
&\ll &WV_0^{1+\ep}  (W^{\ell/2} V^{-2} + 1).\nonumber
\end{eqnarray}
(This improves upon (\ref{NV_weak}) by effectively replacing
$V_0^{-2}$ by $V^{-2}$; observe that up to constant factors, $U$ is the same size as $W^{\ell/2}$.) 
Summing over dyadic regions then shows 
\begin{eqnarray*}
 \sum_{\bstack{V_0 \leq V/2}{\text{dyadic}}} N(Z,X;V_0) &\ll&
 W^{1+\ell/2}V^{-1+\ep} + WV^{1+\ep} \\ [-16pt] 
 	& \ll & X^\ep \{ Z^{2}X^{1/2} + Z^{\ell+2} X^{-1/2} \},
	\end{eqnarray*}
which proves Proposition \ref{prop_Sg_avg}.

\section{Proof of Proposition \ref{prop_Sg_moment}}\label{sec_T1T2}
We define a quantity $R_\ell(d;Z)$ according to the parameters $U,V,W$
given in (\ref{param}) as follows:
set $R_\ell(d;Z)=0$ if $d$ is not square-free, and for $d$ square-free let
$R_\ell(d;Z)$ be the number of triples $(w,u,v)\in\mathbb{N}^3$ satisfying
\[ W\le w<4W,\;\;\;u\le U,\;\;\; v\le V,\;\;\;\gcd(w,v)=1,\]
\[w=p_1p_2\mbox{ with } p_1\not=p_2\in [Z,2Z),\]
and
\[4w^\ell=u^2+dv^2.\]
  Recall also the quantity $S_\ell(d;Z)$ defined in Proposition
  \ref{prop_hg_bound}.  
Upon letting $w=p_1p_2$, we observe that (up to signs) any tuple
$p_1,p_2,u,v$ contributing to $S_\ell(d;Z)$ must have $W \leq w < 4W,  1
\leq u \leq U, 1 \leq v \leq V$, so that $S_\ell (d;Z) \ll R_\ell(d;Z)$. 
Thus we may write
\beq\label{SRR}
 \sum_{X \leq d<2X} S_\ell(d;Z)^2  \ll  \sum_{X \leq d < 2X} R_\ell(d;Z) +
 \sum_{ X \leq d<2X} R_\ell(d;Z)(R_\ell(d;Z)-1). 
 \eeq
The advantage of separating the terms in this fashion is that in the
second term on the right hand side we may now count only distinct
tuples $(u,v,w) \neq (u',v',w')$ in $R_\ell(d;Z)$. 

We note that for $X^{\frac{1}{2\ell}} \le Z \le X$ the first term on the right
hand side of (\ref{SRR}) satisfies  
\beq\label{R_avg}
  \sum_{X \leq d < 2X} R_\ell(d;Z) \ll \sum_{\bstack{V_0 \leq
      V/2}{\text{dyadic}}} N(Z,X;V_0) \ll X^\ep \{ Z^2X^{1/2} +
  Z^{\ell+2} X^{-1/2} \}, 
  \eeq
by Proposition \ref{prop_Sg_avg}. 
The main remaining task is to treat 
\[ T_\ell = T_\ell(Z;X) := \sum_{ X \leq d<2X} R_\ell(d;Z)(R_\ell(d;Z)-1).
 \]
We will prove:
\begin{prop}\label{prop_T1T2}
For $X^{\frac{1}{2\ell}} \le Z \le X$,
\beq \label{T2_X}
T_\ell \ll  Z^{2\ell+4}X^{\ep-1}.
\eeq
Moreover when $\ell =3$ and $X^{\frac{1}{6}}\le Z\le X$ we have
\beq \label{T2_X'}
T_3 \ll  X^{\ep}(Z^7X^{-1/2}+Z^{12}X^{-3/2}).
\eeq
\end{prop}
Combining (\ref{R_avg}) and (\ref{T2_X}), we see that
\[  \sum_{X \leq d<2X} S_\ell(d;Z)^2  \ll X^\ep( Z^2X^{1/2} + Z^{\ell+2}
  X^{-1/2}+Z^{2\ell+4}X^{-1}).\]
Note that 
\[Z^{\ell+2}X^{-1/2}\le Z^{2\ell+4}X^{-1}\]
for $Z\ge X^{1/(2\ell)}$, so that under this assumption 
\[  \sum_{X \leq d<2X} S_\ell(d;Z)^2  \ll X^\ep( Z^2X^{1/2}+Z^{2\ell+4}X^{-1}).\]
This suffices for Proposition \ref{prop_Sg_moment} for $\ell \geq 5$. 
For $\ell =3$ we improve on this; from (\ref{R_avg}) and (\ref{T2_X'}) we obtain
\[  \sum_{X \leq d<2X} S_3(d;Z)^2  \ll X^\ep( Z^2X^{1/2} +Z^5X^{-1/2}+  
Z^7 X^{-1/2}+Z^{12}X^{-3/2}).\]
 However 
\[Z^5X^{-1/2}\le
Z^7X^{-1/2}=\{Z^2X^{1/2}\}^{1/2}\{Z^{12}X^{-3/2}\}^{1/2}
\le Z^2X^{1/2}+Z^{12}X^{-3/2},\]
whence the case $\ell =3$ of
Proposition \ref{prop_Sg_moment} also follows.

\subsection{A first bound for $T_\ell$}\label{sec_T1T2_first}
We now prove (\ref{T2_X}).
We recall the parameters $U,V,W$ of (\ref{param}) and note that $T_\ell$
is at most the number of 6-tuples $(w_1,w_2,u_1,u_2,v_1,v_2)$ in the
ranges 
\[  W \leq w_1,w_2  < 4W,\;\;\; 1 \leq  u_1,u_2 \leq U,\;\;\;
 1 \leq v_1,v_2 \leq V\]
 that satisfy the conditions
 \beq
  (u_1,v_1,w_1) \neq (u_2,v_2,w_2), \label{T2_0} 
\eeq
\beq
\gcd(w_1,v_1) =\gcd(w_2,v_2)=1, \label{T2_1}
\eeq
\beq
 v_1^2 |(4w_1^\ell - u_1^2) \text{ and } v_2^2|(4w_2^\ell - u_2^2), \label{T2_2}
\eeq
\beq  v_1^2(4w_2^\ell - u_2^2) = v_2^2(4w_1^\ell - u_1^2) \neq 0 \label{T2_3}.
  \eeq
  We will obtain a first upper bound for $T_\ell$ by following the
  approach of \cite{Sou00}, ignoring the divisibility conditions
  (\ref{T2_2}); note that we are also now ignoring the fact that each
  of $w_1,w_2$ is a product of two distinct primes. 
We claim that for tuples satisfying the above conditions,
\beq\label{vvww}
v_1^2w_2^\ell-v_2^2w_1^\ell\neq 0. 
\eeq 
To prove this we recall that
$\gcd(w_i,v_i)=1$ for $i=1,2$, whence $v_1^2w_2^\ell=v_2^2w_1^\ell$
would imply that $v_1=v_2$ and $w_1=w_2$, and
 hence $u_1=u_2$. This would then contradict (\ref{T2_0}).

We now observe that once $v_1,v_2,w_1,w_2$ are fixed then $u_1,u_2$ are fixed up to $X^\ep$ choices. For indeed,   fixing $v_1,v_2,w_1,w_2$ in (\ref{T2_3}) gives
\beq\label{key_id}
 4(v_2^2 w_1^\ell - v_1^2w_2^\ell)  = (v_2u_1 - v_1u_2)(v_2u_1 + v_1u_2).
 \eeq
 The left-hand side is a nonzero integer by (\ref{vvww}), 
so that $u_1,u_2$ are fixed up to $X^\ep$ choices.
Thus we obtain  
\[ T_\ell \ll W^2 V^2 X^\ep \ll Z^{2\ell+4}X^{-1+\ep},\]
which is the bound given in (\ref{T2_X}).  

\subsection{A second bound for $T_\ell$}
We may obtain the alternative upper bound (\ref{T2_X'}) for $T_\ell$
by following the method of \cite{HB07a}, but with the
addition of certain technical considerations because in the present
case the variables $v_i$ are not restricted to be primes.  
Although it is easy enough to do this for general odd prime $\ell$ 
we shall confine our attention to $\ell =3$, since this is the only case 
we shall use.

First we consider the contribution to $T_3$ arising from the case in
which $\gcd(w_1,w_2)\not=1$.  We write $T^0_3$ for the number of
6-tuples of this type.  Since each of $w_1$ and $w_2$ is a
product of two primes in the interval $[Z,2Z)$ this can happen only
  when there is at least one prime $p\in[Z,2Z)$ dividing both of $w_1$
    and $w_2$.  The number of possible pairs $w_1,w_2$ is thus
    $O(Z^3)$.  We now follow the argument of  Section
    \ref{sec_T1T2_first}.  There are $O(V^2)$ pairs $v_1,v_2$, and the
    factorization (\ref{key_id}) shows that there are $O(X^{\ep})$
    possibilities for $u_1,u_2$ once $w_1,w_2,v_1,v_2$ are fixed.  It
    follows that
\[T^0_3\ll Z^3V^2X^{\ep}.\]

From now on we assume that $\gcd(w_1,w_2)=1$. For each integer $1 \leq \del \leq V$, 
we will let $T_3(\del)$
denote the contribution to $T_3$ from triples $(u_1,v_1,w_1)$ and
$(u_2,v_2,w_2)$ with $w_1,w_2$ coprime, such that $\gcd(v_1,v_2)=\del$. We
will prove: 

\begin{prop}\label{prop_T2del}
For each integer $1 \leq \del \leq V$,
\[ T_3(\del) \ll X^\ep(W^2 V^{2/3} \del^{-2/3}+ V^3U\del^{-3} 
+ WV\del^{-1}).\] 
\end{prop}
From this we conclude that
\begin{eqnarray*}
 T_3  &\ll& T_3^0 +\sum_{\del=1}^{V} T_3(\del)\\
 & \ll & X^\ep \{Z^3V^2+\sum_{\del=1}^V (W^2 V^{2/3} \del^{-2/3}
 +V^3 U\del^{-3} + WV\del^{-1})\}\\ 
 & \ll & X^\ep (Z^3V^2+VW^2 + V^3 U  +WV ) \\
  & \ll & X^\ep ( Z^3V^2+ VW^2+V^3 U ),
 \end{eqnarray*}
 since clearly $WV \ll VW^2$. 
  Upon recalling the parameter definitions (\ref{param}) this shows that 
\[   T_3 \ll X^\ep\{Z^9X^{-1}+ Z^{7}X^{-1/2}+ Z^{12}X^{-3/2}\}.\]
This provides the second bound for $T_3$ given in 
Proposition \ref{prop_T1T2}, since $Z\ge X^{1/6}$.

To prove Proposition \ref{prop_T2del}, we fix $\del$ and write $v_i  =
\del y_i$ for $i=1,2$ so that $\gcd(y_1,y_2)=1$.  
We first isolate solutions $(u_1,v_1,w_1)$ and $(u_2,v_2,w_2)$ that
contribute to $T_3(\del)$  such that $y_1,y_2$
satisfy a relation
\beq\label{v1v2_ell}
y_1^2 \mu_2^3 =y_2^2 \mu_1^3
\eeq
for some integers $\mu_1,\mu_2$.
Given a relation of the form (\ref{v1v2_ell}), we may divide both sides
by $\gcd(\mu_1,\mu_2)^3$ to obtain an equivalent relation  
\[y_1^2 \lambda_2^3 = y_2^2 \lambda_1^3
\]
in which $(\lambda_1,\lambda_2)=1$ and $(y_1,y_2)=1$.
This implies that for each $i=1,2$,
\beq\label{v_shape}
 y_i^2 =  \lambda_i^3.
 \eeq
 This implies that $y_i$ is itself a perfect cube, say $y_i = s_i^3$. 
We recall from (\ref{key_id}) that once $v_1,v_2, w_1,w_2$ are fixed, $u_1,u_2$ are fixed up to $X^\ep$ choices. 
Thus we count how many $v_1,v_2 \leq V$ with $\gcd(v_1,v_2)=\del$ are
of the type (\ref{v_shape}) by noting that there are at most $O((V\del^{-1})^{1/3})$
choices for each $s_i$. We bound the number of choices for
$w_1,w_2$ trivially by $O(W^2)$, and conclude that the contribution to $T_3(\del)$ of solutions for
which a relation of the form (\ref{v1v2_ell}) holds is at most  
\beq\label{Con1}
\ll W^2  V^{2/3} \del^{-2/3}X^\ep .
 \eeq

We now proceed to count the remaining contribution to $T_3(\del)$; we
may assume from now on that no relation of the form (\ref{v1v2_ell})
holds for $y_1$ and $y_2$. 
Define 
\beq\label{kdfn}
 k = y_2u_1 + y_1u_2.
 \eeq
Note that if $\del, w_1,w_2, y_1,y_2$ and $k$ are fixed, then
$u_1,u_2$ are fixed uniquely by (\ref{key_id}). Thus we will count the
number of solutions $w_1,w_2$ contributing to $T_3(\del)$ for each
fixed $y_1,y_2,k$. 
  
Recalling the definition of $y_1,y_2$ we see that the condition (\ref{T2_3}) now becomes 
\[ y_1^2(4w_2^\ell - u_2^2) = y_2^2(4w_1^\ell - u_1^2) \neq 0,\]
and since $\gcd(y_1,y_2)=1$, this implies a system of congruences
\begin{eqnarray}
4y_2^2 w_1^3 &\con& k^2 \modd{y_1} \label{conv1}\\
4y_1^2 w_2^3 &\con& k^2 \modd{y_2} \label{conv2}\\
4y_2^2 w_1^3 &\con& 4y_1^2 w_2^3 \modd{k}. \label{conk}
\end{eqnarray}
We first reduce this to a similar system of congruences with 
square-free moduli. For  $i=1,2$ let $q_i$ 
denote the odd square-free kernel of $y_i$, that is 
\[ q_i = \prod_{\bstack{p|y_i}{p >2}} p.\]
 The congruence (\ref{conv1}) implies that $4y_2^2 w_1^3 \con k^2 \modd{q_1}$. Since $(4y_2,q_1)=1$ this congruence may be re-written as $w_1^3 \con a_1 \modd{q_1}$ for some constant $a_1$ determined by $y_2$ and $k$. A similar observation applies to (\ref{conv2}).
Next, we define 
\[ r = \prod_{\bstack{p|k}{p >2}}p\]
to be the odd square-free kernel of $k$, and deduce from (\ref{conk}) an analogous congruence modulo $r$. It follows that any solutions $w_1,w_2$ of the system (\ref{conv1})-(\ref{conk}) must satisfy the congruences 
\begin{eqnarray}
w_1^3 &\con& a_1 \modd{q_1} \label{conv1_y}\\
 w_2^3 &\con& a_2 \modd{q_2} \label{conv2_y}\\
y_2^2 w_1^3 &\con& y_1^2 w_2^3 \modd{r} \label{conk_y}
\end{eqnarray}
for some constant $a_1$ determined by $y_2, k \modd{q_1}$ and some constant $a_2 $ determined by $y_1,k \modd{q_2}$.

Certainly $(q_1,q_2)=1$. In addition, we note that $(y_1,r) =1$ and $(y_2,r)=1$. For indeed, if
some odd prime $p$ satisfies $p\mid k$ and $p\mid y_1$ then by (\ref{kdfn})
it follows that $p\mid u_1$, since by construction
$(y_1,y_2)=1$. However, by
the condition $v_1^2 \mid (4w_1^3 - u_1^2)$,  this would imply that
$p|w_1$, which contradicts the fact that $(v_1,w_1)=1$. The fact that
$(y_2,r)=1$ may be shown similarly.  
As a consequence of these observations, 
\beq\label{kvu_star}
(q_1,q_2)=1, \qquad (q_1,r)=1, \qquad (q_2,r)=1.
\eeq
The next step is to note that the conditions (\ref{conv1_y})-(\ref{conk_y}) may be interpreted as lattice conditions.

\begin{lemma}\label{lemma_lattice}
The congruence (\ref{conv1_y}) requires that $w_1$ lies in one of at
most $3^{\om(q_1)}$ residue classes modulo $q_1$, and similarly
(\ref{conv2_y}) requires that $w_2$ lies in one of at most
$3^{\om(q_2)}$ residue classes modulo $q_2$.  

Furthermore, there exists a collection of at most $3^{\om(r)}$
lattices $\Lambda_i \subset \Z^2$ of determinant $r$,
such that any coprime pair $(w_1,w_2)$
satisfying (\ref{conk_y}) must lie in $\Lambda_i$ for
some $i$.  Conversely any pair $(w_1,w_2)$ in any of the lattices
$\Lambda_i$ will satisfy (\ref{conk_y}). 
\end{lemma}

To prove this, we first consider the congruence (\ref{conv1_y}). Fix a prime divisor
$p |q_1$; then $w_1$ can only be a solution to (\ref{conv1_y})  if 
\beq\label{prime}
w_1^3 \con a_1 \modd{p}.
\eeq
There are at most $3$ residue classes modulo
$p$ in which a solution $w_1$ to (\ref{prime}) may lie. We may
conclude that $w_1$ lies in one of at most $3^{\om(q_1)}$ residue
classes modulo $q_1$. A similar argument applies to (\ref{conv2_y}),
establishing that $w_2$ may lie in at most $3^{\om(q_2)}$ residue
classes modulo $q_2$.  

We now turn to (\ref{conk_y}). Since $(y_1,r)=1$ and $(w_1,w_2)=1$
we must have $(w_1,r)=1$. (Indeed, otherwise, if we suppose $p$ is a prime factor of both $w_1$ and $r$, we would 
see in (\ref{conk_y}) that $p | y_1^2w_2^3$, but since $(w_1,w_2)=1$ we cannot have $p | w_2$, so we would conclude $p | y_1$. This would in turn contradict that fact we previously proved that $(y_1,r)=1$.) Using the fact that $(w_1,r)=1$ we see that (\ref{conk_y}) implies
\beq\label{Nw}
w^3\con a\modd{r},
\eeq
where $w\con w_2w_1^{-1}\modd{r}$ and $a\con
(y_2y_1^{-1})^2\modd{r}$ is coprime to $r$.
Now, just as with our analysis of (\ref{conv1_y}), we see that there
is a collection of at most $3^{\omega(r)}$ 
residue classes $w\con b_i\modd{r}$
in which $w$ must lie. This leads to a corresponding collection of
lattice conditions $w_2\con b_iw_1\modd{r}$ which, taken together,
are equivalent to (\ref{Nw}). Finally we note that the resulting lattice of pairs $(w_1,w_2)$ has a basis $\{(1,b_i),(0,r)\}$, so that its determinant is just $r$. This completes the proof of the lemma.

\subsection{Counting lattice points}\label{sec_lattice}

Since $q_1,q_2,r$ are coprime in pairs, we may conclude from
Lemma \ref{lemma_lattice} that $(w_1,w_2)$ must lie in one of  at most
$3^{\om(q_1)+\om(q_2) + \om(r)}$ lattice cosets 
of the form $(c_1,c_2) + \Lambda$,
where $\Lambda$ is a lattice with
$\det(\Lambda) = q_1q_2r$. We note that the total number
of lattices is $\ll X^\ep$, since under 
the assumption $Z \le X$, we
have $v_i \leq V \ll X^{5/2}$ and $k \leq 2UV \ll X^{11/2}$.  We now fix one of these lattices, which we will denote by $\Lambda$, and its corresponding shift $(c_1,c_2)$. Note that
we may choose $(c_1,c_2)$ such that $W \leq c_i < 4W$ for
$i=1,2$, since otherwise $w_1,w_2$ would lie outside the desired range
$W \leq  w_1,w_2 < 4W$.  
We now write $(z_1,z_2) = (w_1,w_2) - (c_1,c_2)$, and proceed to count
the number of  
\[ (z_1,z_2) \in \Lambda, \quad |z_i| < 3W.\]
 Let $\lam_1 \leq \lam_2$ be the successive minima of $\Lambda$, so
 that the standard Minkowski inequalities show that $\det(\Lambda) \ll
 \lam_1 \lam_2 \ll \det(\Lambda)$ (see for example equation (5) of Davenport \cite{Dav58}).
We note that in our particular case, 
\beq\label{lam1}
 \lam_1 \ll \sqrt{\det(\Lambda)} \ll \sqrt{q_1q_2 r} \ll
 V^{3/2} U^{1/2} \del^{-3/2}. 
 \eeq
Here we have used the fact that 
$q_i \leq y_i \leq V \del^{-1}$ for $i=1,2$ and
hence $r \leq k \ll UV\del^{-1}$.  
 Moreover, by Lemma 1 of Davenport \cite{Dav58}, the number of lattice
 points in $\Lam$ with $|(z_1,z_2)| \leq x$  is $ \ll (1+
 x/\lam_1)(1+x/\lam_2).$ 
Thus the number of allowable $z_1,z_2$ in our case is 
\begin{eqnarray*}
& \ll & (1+W/\lam_1)(1+W/\lam_2) \nonumber \\
& \ll & 1 + W^2/ \det (\Lambda) + W/\lam_1 \nonumber \\
& \ll & 1+ W^2/(q_1q_2 r) + W/ \lam_1. 
\end{eqnarray*}
Thus we have 
\beq\label{Tsum}
T_3(\del)  \ll X^\ep \sum_{y_1,y_2,k} (1+
\frac{W^2}{q_1q_2r} + \frac{W}{\lam_1}), 
\eeq
where we recall that $q_i$ is the odd square-free kernel of $y_i$ and 
for each triple $y_1,y_2,k$ we take $\lam_1$ to be the
smallest value from all the corresponding lattices $\Lambda$.  
Recall that $y_1, y_2 \leq V\del^{-1}$ and $k \leq 2 UV\del^{-1}$. 
Then we see that the contribution of the first term in (\ref{Tsum}) to
$T_3(\del)$ is at most  
\beq\label{Con2}
\ll X^\ep V^3 U\del^{-3}.
\eeq
The contribution to
$T_3(\del)$ from the second term in (\ref{Tsum}) is  
\beq\label{T3_sum_3terms}
\ll X^{\ep}W^2\left(\sum_{y_1\le V\del^{-1}}\frac{1}{q_1}\right)
\left(\sum_{y_2\le V\del^{-1}}\frac{1}{q_2}\right)
\left(\sum_{k\le 2UV\del^{-1}}\frac{1}{r}\right).
\eeq
To bound each internal sum we apply the following minor
 variant of Lemma 1 of \cite{HB07a}: 
\begin{lemma}\label{lemma_kappa}
Given an integer $k$, let $k^*$ denote its odd square-free kernel.
For any fixed integer $\kappa \leq K$, 
\[ \#\{ k \leq  K: k^*= \kappa \} \ll K^\ep.\]
\end{lemma}
We defer the proof of this lemma until Section \ref{sec_lemmas}, and merely apply it now to (\ref{T3_sum_3terms}); for example the first sum is bounded by
\[\sum_{y_1\le V\del^{-1}} \frac{1}{q_1}\le \sum_{\nu\leq V\del^{-1}} \frac{1}{\nu} \#\{ v \leq
  V\del^{-1} : v^* = \nu \}\ll V^{\ep}
\sum_{\nu\leq V\del^{-1}} \frac{1}{\nu} \ll V^{\ep}.\]
 One may handle the second and third sums in (\ref{T3_sum_3terms})
similarly, and deduce that the second term in (\ref{Tsum}) is  
$O(W^2X^{\ep})$ overall.  Since $W^2\le W^2V^{2/3}\del^{-2/3}$ for
$\del\le V$ we see that this 
error is dominated by (\ref{Con1}).

Finally, the contribution, say $T_3'(\del)$, of the third term in
(\ref{Tsum}) may be bounded by following the same argument as in
\cite{HB07a}, which we sketch for completeness. For each triple
$y_1,y_2,k$, let $\Lambda$ be the lattice to which $\lam_1$
corresponds, and let $(\mu_1,\mu_2)$ be the shortest non-zero vector
in $\Lambda$, so that $\lam_1$ is the length of $(\mu_1,\mu_2)$. Then  
\[ T_3'(\del) \ll X^\ep W \sum_{\mu_1,\mu_2} \frac{\#
  \{y_1,y_2,k\}}{\sqrt{|\mu_1|^2 + |\mu_2|^2}},\] 
where we count the number of $y_1,y_2,k$ that generate a lattice in which
$(\mu_1,\mu_2)$ is a vector of minimal length. We note by (\ref{lam1}) that  
\beq\label{mu_size_0}
\mu_1,\mu_2 \ll V^{3/2} U^{1/2} \del^{-3/2}.
\eeq
Since $(\mu_1,\mu_2)$ lies in the lattice $\Lambda$, then by construction 
\beq\label{vv}
q_1\mid \mu_1,  \quad q_2 \mid \mu_2
\eeq
and 
\beq\label{kuu}
 r \mid (y_2^2 \mu_1^3 - y_1^2 \mu_2^3),
 \eeq
as described in Lemma \ref{lemma_lattice}.

We first consider the case where both $\mu_1,\mu_2$ are nonzero. By
(\ref{vv}), once $\mu_1,\mu_2$ are fixed, they determine at most
$X^\ep$ values of $q_1, q_2$ and hence at most $X^\ep$ values for
$y_1,y_2$ by Lemma \ref{lemma_kappa}. If $y_2^2 \mu_1^3 -
y_1^2\mu_2^3$ is nonzero, then it determines at most $X^\ep$ 
values for $r$ by (\ref{kuu}) and hence at most $X^\ep$ values for
$k$.  
On the other hand, if 
\beq\label{vuvu}
y_2^2 \mu_1^3 =y_1^2 \mu_2^3,
\eeq
then $y_1,y_2$ would satisfy a relation of the form (\ref{v1v2_ell});
pairs $y_1,y_2$ of this type have already been treated, and are
excluded from the contribution we are currently calculating.  
We therefore see that the contribution to $T_3'(\del)$ from
$\mu_1,\mu_2$ both nonzero is 
\[ T_3'(\del) \ll  X^{4\ep} W \sum_{\mu_1,\mu_2}
\frac{1}{\sqrt{|\mu_1|^2 + |\mu_2|^2}} .\] 

To bound the sum, we begin by focusing on a fixed dyadic range 
\[\frac{1}{2}B < \sqrt{|\mu_1|^2 + |\mu_2|^2} \leq B,
\]
for any appropriate $B \geq 1$;
we note that the restriction (\ref{mu_size_0}) implies that $B \ll
V^{3/2}U^{1/2}\del^{-3/2}.$  There are $O(B^2)$ pairs $\mu_1,\mu_2$, each of
which contribute $O(B^{-1})$ to the sum. 
 Summing over dyadic $B\ll
V^{3/2}U^{1/2}\del^{-3/2}$ therefore produces a 
total contribution of $\ll
X^{\ep}WV^{3/2}U^{1/2}\del^{-3/2}$ to $T_3'(\del)$.

On the other hand if $\mu_1$ vanishes, then there are $V\del^{-1}$ choices
for $y_1$ and $O(X^{2\ep})$ choices for $q_2, r$, hence
$O(X^{4\ep})$ choices for $y_2,k$. (In particular, (\ref{vuvu})
cannot occur, since it would force $\mu_1 = \mu_2=0$.) Thus the
contribution from these terms to $T_3'(\del)$ is  
\[ \ll X^{5\ep}  V W\del^{-1} \sum_{\mu_2 \ll V^{3/2} U^{1/2} \del^{-3/2}}
\frac{1}{|\mu_2|} \ll X^{6\ep} VW \del^{-1}.\] 
The case where $\mu_2$ vanishes may be treated by an analogous
argument. We may conclude that  
\[ T_3'(\del) \ll  X^\ep (W V^{3/2} U^{1/2} \del^{-3/2} + VW \del^{-1}) .\]
Combining this with the contributions (\ref{Con1}) and (\ref{Con2})
shows that  
\[ T_3(\del) \ll X^\ep(W^2 V^{2/3}\del^{-2/3} + V^3 U\del^{-3}
+ W V^{3/2} U^{1/2} \del^{-3/2} + VW\del^{-1}).\]
Since
\begin{eqnarray*}
W V^{3/2} U^{1/2} \del^{-3/2} & = & \left\{W^2\right\}^{1/2}
\left\{V^3U\del^{-3}\right\}^{1/2}\\
&\le&\left\{W^2V^{2/3}\del^{-2/3}\right\}^{1/2}
\left\{V^3U\del^{-3}\right\}^{1/2}
\end{eqnarray*}
for $\del\le V$, the third term above is dominated by the first two, so
that Proposition \ref{prop_T2del} follows.

\subsection{Proof of Lemma \ref{lemma_kappa}}\label{sec_lemmas}
We now prove Lemma \ref{lemma_kappa}, in the following more general
form. Given any finite set $\Pcal$ of primes (possibly empty), let  
\[ k(\Pcal) = \prod_{\bstack{p |k}{p \not\in \Pcal}} p.\] 
Consider the set $\{ k \leq K : k(\Pcal) = \kappa\}$ for a fixed
positive integer $\kappa$. The set is empty unless $\kappa \leq K$ is
square-free and satisfies $(\kappa, \prod_{p \in \Pcal} p) =1$, which
we now assume.  Then  for any $\eta>0$,  
\begin{eqnarray*}
\#\{ k \leq K : k(\Pcal)= \kappa \}
	& \leq & \sum_{\bstack{k=1}{k(\Pcal)= \kappa} }^K \left(
          \frac{K}{k} \right)^\eta \\ 
	& \leq & K^\eta  \sum_{\bstack{k=1}{k(\Pcal) = \kappa}}^\infty k^{-\eta} \\
	& = & K^\eta  \prod_{p \in \Pcal} \left( \sum_{e=0}^\infty
          p^{-e \eta} \right) \prod_{p |\kappa} \left(
          \sum_{e=1}^\infty p^{-e \eta} \right). 
\end{eqnarray*}
Setting $A(\eta) = \sum_{e=0}^\infty 2^{-e\eta}$ we then see that 
\[ \#\{ k \leq K : k(\Pcal) = \kappa \} \leq K^\eta
A(\eta)^{\om(\kappa)+\#\Pcal} \leq K^\eta A(\eta)^{(\#\Pcal +
  1)\om(\kappa)}.  \] 
Upon recalling that $\om(\kappa) \ll (\log 3 \kappa)(\log \log 3
\kappa)^{-1}$ and $\kappa \leq K$ we may conclude that  
\[ \#\{ k \leq K : k(\Pcal)= \kappa \} \ll_\eta K^{(\#\Pcal+ 2)\eta} \]
for any $\eta>0$, which proves Lemma \ref{lemma_kappa}.

\section{Average of $h_\ell(-d)$}\label{sec_avg_h}
We now turn to applications of the key propositions. We first apply
Proposition \ref{prop_hg_bound} to derive a nontrivial upper bound for
the average of $h_\ell(-d)$. Fix a dyadic region $X \leq d < 2X$ and
assume that $X^{1/(2\ell)} \le Z \le X$. Then Proposition
\ref{prop_hg_bound} implies that 
\begin{eqnarray*}
\lefteqn{\sum_{X\leq d<2X} h_\ell(-d)}\\
& \ll& X^\ep \{X^{1/2}\# E(Z;X) + X^{3/2} Z^{-1} +
         X^{1/2} Z^{-2}\sum_{\bstack{X\leq d < 2X}{d \not\in E(Z;X)}}  S_\ell(d;Z) \} .
\end{eqnarray*} 
We apply the upper bound (\ref{E_small}) to the exceptional set
$E(Z;X)$ and Proposition \ref{prop_Sg_avg} to the average of
$S_\ell(d;Z)$ to conclude that 
\[\sum_{X \leq d<2X} h_\ell(-d) 
	\ll  X^{\ep}\{ X^{3/2}Z^{-1} + X +  Z^\ell \}.
	\]
It is optimal to choose $Z=X^{\frac{3}{2\ell+2}}$, resulting in
\[\sum_{X\leq d<2X} h_\ell(-d)  \ll X^{\frac{3}{2} - \frac{3}{2\ell+2} + \ep}.\]
 Summing over $O(\log X)$ dyadic intervals to cover the full range
 $0<d< X$ then yields the result of Theorem \ref{thm_avg}.

\section{Higher moments of $h_\ell(-d)$}\label{sec_higher}
\label{sec_second}  
We now consider higher moments.
For any odd prime $\ell$, define for any real $H \geq 1$ the set
\[ A_\ell(H;X) = \{ X\le d<2X : h_\ell(-d) > H \},\]
with corresponding counting function
\[ N_\ell(H;X) = \# A_\ell(H;X).\]
We also define for any $\frac{1}{4}X^{\frac{1}{2\ell}} \leq Z \leq X$ the set
\[ A_\ell^0(H,Z;X) = \{ X\leq d<2X : h_\ell(-d) > H \} \setminus E(Z;X),\]
where $E(Z;X)$ is as usual the exceptional set provided by Proposition \ref{prop_hg_bound}.
We define the corresponding counting function
\[ N_\ell^0(H,Z;X) = \# A_\ell^0(H,Z;X).\]
We note that for any fixed choice of $Z$ in the above range,
\beq\label{NEN}
N_\ell(H;X) \le  \#E(Z;X) +N_\ell^0(H,Z;X) \ll X^\ep  + N_\ell^0(H,Z;X). 
\eeq

\subsection{The case $\ell =3$}
Restricting to the case $\ell =3$, we see that (\ref{DavHeil_bound}) implies that 
\beq\label{N3H}
N_3(H;X) \ll XH^{-1}.
\eeq
We also note that $A_3(H;X)$ is empty by (\ref{EllVen}) unless $H
\le X^{1/3+\ep}$ for some small $\ep>0$. 
In general we have:
\begin{prop}\label{prop_NH3}
For $1 \leq H \leq X^{1/3+\ep}$,
\[ N_3(H;X) \ll X^\ep (X^{1/2}+ X^{7/2}H^{-10}).\]
\end{prop}

To prove this we consider $A_3^0(H,Z;X)$ with the choice $Z= X^{1/2 + 2\ep} H^{-1}$; note in
particular $Z \ge X^{1/6}$ when $H\le X^{1/3+\ep}$. 
Moreover we will have
\[ h_3(-d) > H \gg d^{1/2 + \ep} Z^{-1}\]
for all $d$ in $A_3^0(H,Z;X),$  whence
Proposition \ref{prop_hg_bound} shows that
\[ h_3(-d) \ll d^{1/2+\ep}Z^{-2}  S_3(d;Z).\]
We therefore have
\[ S_3(d;Z) \gg d^{-1/2-\ep} Z^{2} h_3(-d) \gg X^{-1/2 -\ep}Z^2
h_3(-d) >  X^{-1/2 -\ep}Z^2 H, \] 
for all $d \in A_3^0(H,Z;X)$.  This leads to the bound
\[ N_3^0(H;X)\left(X^{-1/2 -\ep}Z^2 H\right)^2 \ll \sum_{d \in  A_3^0(H,Z;X)}S_3(d;Z)^2
 \ll \sum_{X \leq d <2X} S_3(d;Z)^2.  
 \]
We can now apply the case $\ell =3$ of 
Proposition \ref{prop_Sg_moment} to obtain
\[ N_3^0(H,Z;X)\left(X^{-1/2 -\ep}Z^2 H\right)^2 
\ll X^{\ep}\{Z^2X^{1/2}+Z^{12}X^{-3/2}\},\]
so that 
\[ N_3^0(H,Z;X) \ll X^{3\ep}H^{-2}\{Z^{-2}X^{3/2}+Z^8X^{-1/2}\}
\ll X^{19\ep}\{X^{1/2}+X^{7/2}H^{-10}\}\]
in view of our choice of $Z$. 
This is sufficient for Proposition \ref{prop_NH3}, by (\ref{NEN}).

We may now derive Theorem \ref{thm_square} from Proposition \ref{prop_NH3}.
It will suffice to consider a dyadic range $X \leq d < 2X$. Then
\begin{eqnarray*}
 \sum_{X\le d<2X} h_3(-d)^k &\ll& \sum_{\bstack{H \le
     X^{1/3+\ep}}{\text{dyadic}}} \sum_{\bstack{X\le d<2X}{H < h_3(-d)
     \leq 2H}} h_3(-d)^k \\ 
	&\leq &\sum_{\bstack{H \le X^{1/3+\ep}}{\text{dyadic}}} N_3(H;X)(2H)^k.
	\end{eqnarray*}
In view of (\ref{N3H}) we have
\[N_3(H;X)(2H)^k\ll XH^{k-1}.\]
On the other hand, Proposition \ref{prop_NH3} yields
\[N_3(H;X)(2H)^k\ll X^\ep (X^{1/2}H^{k}+ X^{7/2}H^{k-10}).\]
In particular for $k=4$ we deduce that
\begin{eqnarray*}
N_3(H;X)(2H)^4&\ll& X^\ep\min\{XH^3\,,\,X^{1/2}H^4+ X^{7/2}H^{-6}\}\\
&\ll& X^\ep\min\{XH^3\,,\,X^{1/2}H^4\}+ \min\{XH^3\,,\,X^{7/2}H^{-6}\}.
\end{eqnarray*}
For $H\le X^{1/3+\ep}$ the first term is at most 
\[X^{1/2}H^4\le X^{11/6+4\ep}\]
while the second term is at most
\[\{XH^3\}^{2/3}\{X^{7/2}H^{-6}\}^{1/3}=X^{11/6}.\]
It follows that $N_3(H;X)(2H)^4\ll X^{11/6+4\ep}$, whence
\[ \sum_{X\le d<2X} h_3(-d)^4\ll X^{11/6+5\ep}.\]
This suffices for Theorem \ref{thm_square}.  As noted in the
introduction, one can deduce estimates for other moments from the
fourth moment.  The reader may check that a direct application of the
methods of this section to the general moment only reproduces these
consequences of the special case $k=4$.

\subsection{The case $\ell \geq 5$}
We now consider the $k$-th moment of $h_\ell(-d)$ for primes $\ell \geq 5$
and any real $k \geq 1$.  
By Corollary \ref{newcor}
we see that  
\[ N_\ell(H;X) \ll X^{\ep}  \quad \text{if $H \ge X^{\frac{1}{2} - \frac{1}{2\ell}+\ep}.$} 
\]
We also record the trivial bound
\beq\label{N_triv}
N_\ell(H;X) \ll X,
\eeq
valid for all $H$.  In addition, we claim:
\begin{prop}\label{prop_NH}
For any prime $\ell \geq 3$ and $1 \leq H \leq X^{1/2  - 1/(2\ell)+\ep}$,
\[ N_\ell(H;X) \ll X^\ep (XH^{-1} + X^{\ell/2}H^{-(\ell+1)}).\]
\end{prop}
With Proposition \ref{prop_NH} in hand, we will prove:
\begin{prop}\label{prop_NH_h}
For any prime $\ell \geq 5$ and any real number $k \geq 1$,
\[ \sum_{X\le d<2X} h_\ell(-d)^k \ll X^{\sigma + \ep},\]
where 
\[ \sig = \max \{ \sig_1,\sig_2, \sig_3\} \]
and
\begin{eqnarray*}
\sig_1 & = & 1 + k \left( \frac{\ell-2}{2\ell+2} \right),\\
\sig_2 & = & 1 + k\left( \frac{\ell-1}{2\ell} \right) - \left(
  \frac{\ell-1}{2\ell} \right),\\
\sig_3 & = & \frac{k}{2}\,.
 \end{eqnarray*} 
\end{prop}
We note that the maximum is $\sig_1$ in the range $1 \leq k \leq
\frac{\ell^2-1}{2\ell-1}$; it is $\sig_2$ in the range  
$\frac{\ell^2-1}{2\ell-1} \leq k \leq \ell+1$; and it is $\sig_3$ for $k\geq \ell+1$.
This leads immediately to the
statement of Theorem \ref{thm_k_real}. 
We note that Proposition \ref{prop_NH} does not imply any new results
in the case of $h_3(-d)$.

The proof of Proposition \ref{prop_NH} follows similar lines to that of Proposition \ref{prop_NH3}.  As before we 
set $Z= X^{1/2 + 2\ep} H^{-1}$,
so that $Z \geq X^{1/(2\ell)}$ for $H \leq X^{1/2 - 1/(2\ell) + \ep}$.
We deduce that
\[ S_\ell(d;Z) \gg d^{-1/2-\ep} Z^{2} h_\ell(-d) \gg X^{-1/2 -\ep}Z^2
h_\ell(-d) >  X^{-1/2 -\ep}Z^2 H, \] 
again under the assumption that $d \in A_\ell^0(H,Z;X)$.
As a result,
\[ N_\ell^0(H,Z;X) X^{-1/2 -\ep}Z^2 H \ll \sum_{d \in  A_\ell^0(H,Z;Z)}S_\ell(d;Z)
 \ll \sum_{X \leq d <2X} S_\ell(d;Z).  \]
Upon applying Proposition \ref{prop_Sg_avg} we obtain
\[ N_\ell^0(H,Z;X)X^{-1/2-\ep} Z^2 H \ll X^\ep( Z^2 X^{1/2} + Z^{\ell+2}X^{-1/2}),\]
so that 
\[ N_\ell^0(H,Z;X) \ll X^{2\ep}( X H^{-1} + Z^\ell H^{-1} ) \ll X^{(2+2\ell)\ep}
(XH^{-1} + X^{\ell/2}H^{-(\ell+1)}),\] 
upon recalling the choice of $Z$. 
This is sufficient for Proposition \ref{prop_NH}, by (\ref{NEN}).

We  turn finally to Proposition \ref{prop_NH_h}, for which we initially
fix any real number $k \geq 1$. 
We have already observed that $N_\ell(H;X) \ll X^{\ep}$ if 
\[X^{1/2 -
  1/(2\ell)+\ep} \le H \le X^{1/2+\ep},\]
   which shows that for such $H$, 
\beq\label{N_H_small}
 N_\ell(H;X)H^k \ll X^{k/2 +\ep}  .
 \eeq
  Thus we now instead assume that 
\beq\label{H_bound}
H \le X^{1/2 - 1/(2\ell)+\ep}.
\eeq
Then by the trivial bound (\ref{N_triv}) and Proposition \ref{prop_NH} we have
\begin{eqnarray*}
 N_\ell(H;X) H^k &\ll& X^\ep \min \{XH^k, XH^{k-1} + X^{\ell/2} H^{k-\ell-1} \} \\
 & \ll & X^\ep ( XH^{k-1} + \min \{ XH^k, X^{\ell/2} H^{k-\ell-1} \}).
 \end{eqnarray*}
Under (\ref{H_bound}), the first term is $\ll X^{\sig_2}$.
As long as $k \leq \ell+1$, the second term is largest when $XH^k =
X^{\ell/2} H^{k-\ell-1}$, namely when 
\[ H=X^{\frac{\ell-2}{2\ell+2}} = X^{\frac{1}{2} - \frac{3}{2\ell+2}}.\]  
We may conclude that if $k\le \ell+1$ and $H\le X^{1/2 - 1/(2\ell)+\ep}$ then
\[
 N_\ell(H;X)H^k \ll X^\ep(X^{\sig_1} + X^{\sig_2}),
\]
 with the notation of Proposition \ref{prop_NH_h}. On the other hand, 
if $k\ge \ell+1$ then
\[X^{\ell/2}H^{k-\ell-1}\le X^{\ell/2}H^{k-\ell}\le X^{\ell/2}(X^{1/2})^{k-\ell}=X^{k/2}.\]
Thus $ N_\ell(H;X)H^k \ll X^\ep(X^{\sigma_2} + X^{k/2})$ in this case; note that the second term dominates in the range $k \geq \ell+1$. 
To conclude, 
\beq\label{N_H_next}
 N_\ell(H;X)H^k \ll X^\ep(X^{\sig_1} + X^{\sig_2}+X^{\sig_3})
 \eeq
for all $k\ge 1$.

Combining (\ref{N_H_small}) and (\ref{N_H_next}) shows that 
\begin{eqnarray*}
 \sum_{X\le d<2X} h_\ell(-d)^k &\ll& \sum_{\bstack{H \ll
     X^{1/2+\ep}}{\text{dyadic}}} \sum_{\bstack{X\le d<2X}{H < h_\ell(-d)
     \leq 2H}} h_\ell(-d)^k \\ 
	&\leq &\sum_{\bstack{H \ll X^{1/2+\ep}}{\text{dyadic}}} N_\ell(H;X)(2H)^k \\
	& \ll&  X^{\ep}(X^{\sig_1} + X^{\sig_2}+X^{\sig_3}).
	\end{eqnarray*}
We note that $k/2 \leq \max\{\sig_1,\sig_2\}$ in the range $k \leq \ell+1$.
This proves Proposition \ref{prop_NH_h}, and hence Theorem \ref{thm_k_real}.
The reader may verify that a similar computation based on Proposition \ref{prop_Sg_moment} yields no improvements.

\subsection*{Acknowledgements} 
The authors thank Peter Sarnak for asking a question that spurred this line of enquiry.
The first author was supported by EPSRC grant number EP/K021132X/1. The second author was partially supported by NSF DMS-1402121, and thanks the Hausdorff Center for Mathematics for a very pleasant working environment.

\bibliographystyle{alpha}
\bibliography{NoThBibliography}

\begin{thebibliography}{BST13}

\bibitem[BS96]{BruSil96}
A.~Brumer and J.~H. Silverman.
\newblock The number of elliptic curves over {$\bold Q$} with conductor {$N$}.
\newblock {\em Manuscripta Math.}, 91(1):95--102, 1996.

\bibitem[BST13]{BST}
M.~Bhargava, A.~Shankar, and J.~Tsimerman.
\newblock On the {D}avenport-{H}eilbronn theorem and second order terms.
\newblock {\em Invent. Math.}, 193:439--499, 2013.

\bibitem[CL84]{CohLen84}
H.~Cohen and H.~W. Lenstra, Jr.
\newblock Heuristics on class groups of number fields.
\newblock In {\em Number theory, {N}oordwijkerhout 1983 ({N}oordwijkerhout,
  1983)}, volume 1068 of {\em Lecture Notes in Math.}, pages 33--62. Springer,
  Berlin, 1984.

\bibitem[Dav58]{Dav58}
H.~Davenport.
\newblock Indefinite quadratic forms in many variables {II}.
\newblock {\em Proc. London Math. Soc. (3)}, 8:109--126, 1958.

\bibitem[Dav00]{Dav}
H.~Davenport.
\newblock {\em Multiplicative {N}umber {T}heory}.
\newblock Graduate Texts in Mathematics 74, Springer Verlag, 3rd edition, 2000.

\bibitem[DH71]{DavHei71}
H.~Davenport and H.~Heilbronn.
\newblock On the density of discriminants of cubic fields {II}.
\newblock {\em Proc. Roy. Soc. Lond. A.}, 322:405--420, 1971.

\bibitem[Duk98]{Duk98}
W.~Duke.
\newblock Bounds for arithmetic multiplicities.
\newblock In {\em Proceedings of the {I}nternational {C}ongress of
  {M}athematicians, {V}ol. {II} ({B}erlin, 1998)}, number Extra Vol. II, pages
  163--172, 1998.

\bibitem[EV07]{EllVen07}
J.~S. Ellenberg and A.~Venkatesh.
\newblock Reflection principles and bounds for class group torsion.
\newblock {\em Int. Math. Res. Not. IMRN}, (1):Art. ID rnm002, 18, 2007.

\bibitem[HB07]{HB07a}
D.~R. Heath-Brown.
\newblock Quadratic class numbers divisible by 3.
\newblock {\em Funct. Approx. Comment. Math.}, 37(1):203--211, 2007.

\bibitem[Hou10]{Hou13}
R.~Hough.
\newblock Average equidistribution of {H}eegner points associated to the 3-part
  of the class group of imaginary quadratic fields.
\newblock {\em arXiv:1005.1458v2}, 2010.

\bibitem[HV06]{HelVen06}
H.~A. Helfgott and A.~Venkatesh.
\newblock Integral points on elliptic curves and 3-torsion in class groups.
\newblock {\em J. Amer. Math. Soc.}, 19(3):527--550, 2006.

\bibitem[Sch32]{Sch32}
A.~Scholz.
\newblock \"{U}ber die {B}eziehung der {K}lassenzahlen quadratischer
  {K}\"{o}rper.
\newblock {\em J. Reine Angew. Math.}, 166:201--203, 1932.

\bibitem[Sou00]{Sou00}
K.~Soundararajan.
\newblock Divisibility of class numbers of imaginary quadratic fields.
\newblock {\em J. London Math. Soc. (2)}, 61(3):681--690, 2000.

\bibitem[TT13]{TanTho13}
T.~Taniguchi and F.~Thorne.
\newblock The secondary term in the counting function for cubic fields.
\newblock {\em Duke Math. J.}, 162:2451--2508, 2013.

\bibitem[Zha05]{Zha05}
S.-W. Zhang.
\newblock Equidistribution of {CM}-points on quaternion {S}himura varieties.
\newblock {\em Int. Math. Res. Not.}, (59):3657--3689, 2005.

\end{thebibliography}

 \label{endofpaper}

\end{document}